\documentclass[12pt]{article}

\usepackage{amsmath,amsthm,amssymb}
\usepackage{tikz}
\usepackage[colorlinks,
            linkcolor=red,
            anchorcolor=blue,
            citecolor=blue
            ]{hyperref}
\usepackage{pdfsync}
\usepackage[numbers,sort&compress]{natbib}
\allowdisplaybreaks

\theoremstyle{definition}

\theoremstyle{plain}
\newtheorem{theorem}{Theorem}[section]

\newtheorem{remark}{Remark}[section]

\newtheorem{proposition}{Proposition}[section]

\newtheorem{corollary}{Corollary}[section]

\numberwithin{equation}{section}

\begin{document}
	\title{\bf Harmonic functions and comparsion geometry for noncompact 3-manifolds with nonnegative scalar curvature }

 \author{Jiannan Chen$^{1}$\footnote{Email address: chenjiannan@email.ncu.edu.cn}, Haiping Fu$^{1}$\footnote{Email address: mathfu@126.com\,\,\    
 	}\\
 	{\small 1.Department of Mathematics,  Nanchang University,} {\small Nanchang 330031, Jiangxi, China.} 
    \\
}

\date{}

\maketitle
\begin{abstract}
In this paper, we establish several sharp comparsion results for noncompact 3-manifolds with nonnegative scalar curvature and with a compact, connected boundary under a boundary mean curvature condition. Then we apply these gradient estimates to establish a Cohn-Vossen type bound for the average of the scalar curvature. This generalizes the work of Colding-Minicozzi \cite{Cold4} and Oronzio \cite{Oro2}.
\end{abstract}	

	

\section{\bf Introduction and main results}
Harmonic functions play a significant role in geometric analysis. A classical result is Colding's work \cite{Cold1}, which establishes new monotonicity formulas  along the level sets of harmonic functions via gradient estimates of the Green's function in nonparbolic Riemannian manifolds with nonnegative Ricci curvature. Subsequently, using the same monotone quantities, Agostiniani-Fogagnolo-Mazzieri \cite{Fog1} obtained the new Willmore-type geometrical inequalities. Furthermore, in \cite{Ben},  Benatti-Fogagnolo-Mazzieri generalized the monotone quantities in \cite{Cold1} to establish an optimal version of the Minskowski inequality in nonparbolic Riemannian manifolds with nonnegative Ricci curvature and Euclidean volume growth. 

In recent years, the level sets of harmonic functions have been used to investigate 3–manifolds with nonnegative scalar curvature. A pioneering work of Stern \cite{Stern} discovered intriguing analogy between the use of such level sets and the use of stable minimal surfaces instituted by Schoen and Yau \cite{SY1}. Since then, harmonic functions have found a number of applications in the study of asymptotically flat manifolds with nonnegative scalar curvature. In fact, several new proofs of the positive mass theorem were found in \cite{Bray2} and \cite{Ag3} by means of harmonic functions of different types. Later on, Miao \cite{Miao} derived three monotone quantities along the level sets of a harmonic function and obtain several mass-capacity inequalities. Moreover, he also found several new proofs of the positive mass Theorem.  In \cite{Xia}, Xia-Yin-Zhou established general monotone quantities associated to p-capacitary functions whose model space is a spatial Schwarzschild manifold outside a rotationally symmetric sphere via an ODEs analysis. By using an ODEs technique, Mazurowski-Yao \cite{Maz} obtained the similar monotone quantities independently. For more relevant results, refer to \cite{Mazz,Bray3,Hrisch1,Hrisch2,Miao 0,Miao1}.

Ouside the asymptotically flat setting, gradient integral estimates along the Green's function are a powerful tool to investigate 3-manifolds with nonnegative scalar curvature. This comparsion result was first obtained by Munteanu-Wang in \cite{Wang1}, under the additional assumption of asymptotically nonnegative Ricci curvature. Later on, Chan-Chu-Lee-Tsang \cite{Chan} proved the similar results provided the Ricci curvature is bounded from below. Furthermore, Chodosh-Li \cite{Otis} generalized Munteanu-Wang's monotonicity formula to 3-manifolds that do not have nonnegative scalar curvature, using this, they resolved the stable Bernstein problem in $ \mathbb{R}^{4} $. Recently, Colding-Minicozzi \cite{Cold4} established the sharp gradient integral estimates to obtain the sharp area comparsion, under the assumption that the norm of gradient vanishes at infinity. Moreover, they applied these gradient estimates to derive a Cohn-Vossen type bound for the average of the scalar curvature. Shortly after, Oronzio \cite{Oro2} generalized the work of Colding–Minicozzi \cite{Cold4} to 3–manifolds with a compact, connected, minimal boundary and established several sharp comparsion results. 

In this paper, we generalize the work of Colding-Minicozzi \cite{Cold4} and Oronzio \cite{Oro2} to 3-manifolds with nonnegative scalar curvature and with a compact, connected boundary under a boundary mean curvature condition. More precisely, we establish several sharp comparsion results for 3-manifolds with nonnegative scalar curvature and with a compact, connected boundary under the boundary mean curvature condition (\ref{H condition0 0}). Then we apply these gradient estimates to establish a Cohn-Vessen type bound for the average of the scalar curvature. The main novelty is the boundary mean curvature condition (\ref{H condition0 0}), which is inspired by the work of Chen-Fu \cite{Chen} and Oronzio \cite{Oro1,Oro2}. This condition is key to present a contradiction argument to prove the main results.

In order to state our main results, we first introduce several preliminaries.

The spatial Schwarzschild manifold of mass $ \mathfrak{m} $ outside a rotationally symmetric sphere is the 3-dimensional Riemannian manifold with boundary given by the couple 
\begin{equation}\label{Sch}
\Big(\mathbb{R}^{3}\setminus B_{r_{0}}(0), (1+\frac{\mathfrak{m}}{2|x|})^{4}g_{eucl} \Big),
\end{equation}
which will be denoted by $ (M_{m,r_{0}},g_{m}) $.

Let $ (M,g) $ be an one-ended manifold with boundary $ \Sigma $. Let $ u $ be the solution of the Dirichlet problem:
\begin{equation}\label{u}
\begin{cases}
\Delta u=0 &\text{in } M,\\
u=0, &\text{on } \Sigma,\\
u\to 1, &\text{at } \infty.
\end{cases}
\end{equation}
For any $ k\in (0,2] $, we denote the level set of $ u $ 
$$ 
\Sigma_{t}=\{x\in M| u(x)=f(t)\}, 
$$ 
where 
$$ f(t)=1-\frac{\mathfrak{c}_{\Sigma}}{t+\frac{k-1}{k}\mathfrak{c}_{\Sigma}},\,\,\,\,t\in [\frac{\mathfrak{c}_{\Sigma}}{k}, +\infty).
$$
Set
$$
\mathcal{T}=\Big\{t\in [\frac{\mathfrak{c}_{\Sigma}}{k}, +\infty)\,\,\Big|\,\, f(t) \text{ is a regular value of } u \Big\}.
$$
Let $ \mathfrak{c}_{\Sigma}>0 $ be the boundary capacity of $ \Sigma $ in $ (M,g) $, that is 
$$
\mathfrak{c}_{\Sigma}=\frac{1}{4\pi}\int_{M}|\nabla u|^{2}=\frac{1}{4\pi}\int_{\Sigma}|\nabla u|.
$$
We state the main results. 
\begin{theorem}\label{Th1.1}
Let $ (M,g) $ be a 3-dimensional, complete, noncompact and one-ended Riemannian manifold with a compact and connected boundary $ \Sigma $ and with nonnegative scalar curvature. Suppose that $ H_{2}(M,\Sigma)=0 $ and that there exists a smooth solution $ u $ satisfying (\ref{u}) and $ |\nabla u| $ vanishes at infinity. For any $ k\in (0,2] $, if there exists 
\begin{equation}\label{a condi 0}
\alpha\in (-\frac{1}{k\mathfrak{c}_{\Sigma}},\frac{1}{k\mathfrak{c}_{\Sigma}}) \end{equation}
such that 
\begin{equation}\label{H condition0 0}
H\leq \alpha(1-k^{2}\mathfrak{c}_{\Sigma}|\nabla u|)\,\,\,\text{on }\Sigma.
\end{equation}
Then we have the following results:
\\
$ (a) $ (Boundary gradient integral estimates ) There holds
$$
\int_{\Sigma}|\nabla u|^{2}\leq \frac{4\pi}{k^{2}}.
$$
If the equality holds, then $ (M,g) $ is isometric to $ (M_{m,r_{0}},g_{m}) $ of mass $ \mathfrak{m}=2r_{0}(k-1) $ with $ r_{0}=\frac{\mathfrak{c}_{\Sigma}}{k} $.
\\
$ (b) $ For all $ t\in [\frac{\mathfrak{c}_{\Sigma}}{k},+\infty ) $,
$$
\frac{t^{2}}{\mathfrak{c}_{\Sigma}^{2}}\Big(1+\frac{k-1}{kt}\mathfrak{c}_{\Sigma}\Big)^{4}\int_{\Sigma_{t}}|\nabla u|^{2}\leq 4\pi.
$$
If the equality holds for some $ t_{0}\in [\frac{\mathfrak{c}_{\Sigma}}{k}, +\infty) $, then $ (M,g) $ is isometric to $ (M_{m,r_{0}},g_{m}) $ of mass $ \mathfrak{m}=2r_{0}(k-1) $ with $ r_{0}=\frac{\mathfrak{c}_{\Sigma}}{k} $.
\\
$ (c) $ (Area comparsion) For any $ t\in [\frac{\mathfrak{c}_{\Sigma}}{k},+\infty ) $, 
$$
|\Sigma_{t}^{*}|\geq 4\pi t^{2}\Big(1+\frac{k-1}{kt}\mathfrak{c}_{\Sigma}\Big)^{4},
$$
where $ \Sigma_{t}^{*}=\Sigma_{t}\setminus\{|\nabla u|=0\} $. If the equality holds for some $ t_{0}\in [\frac{\mathfrak{c}_{\Sigma}}{k}, +\infty) $, then $ (M,g) $ is isometric to $ (M_{m,r_{0}},g_{m}) $ of mass $ \mathfrak{m}=2r_{0}(k-1) $ with $ r_{0}=\frac{\mathfrak{c}_{\Sigma}}{k} $.
\\
$ (d) $ (Area-capacity inequality) There holds
$$
\sqrt{\frac{|\Sigma|}{16\pi}}\geq \mathfrak{c}_{\Sigma}.
$$
If the equality holds, then $ (M,g) $ is isometric to $ (M_{m,r_{0}},g_{m}) $ of mass $ \mathfrak{m}=2r_{0}(k-1) $ with $ r_{0}=\frac{\mathfrak{c}_{\Sigma}}{k} $.
\\
$ (e) $ (Volume comparsion) For every $ t\in [\frac{\mathfrak{c}_{\Sigma}}{k},+\infty ) $, 
$$
Vol\{u\leq f(t)\}\geq \int_{\frac{\mathfrak{c}_{\Sigma}}{k}}^{t}4\pi s^{2}\Big(1+\frac{k-1}{ks}\mathfrak{c}_{\Sigma}\Big)^{6}ds.
$$
If the equality holds for some $ t_{0}\in [\frac{\mathfrak{c}_{\Sigma}}{k}, +\infty) $, then $ (M,g) $ is isometric to $ (M_{m,r_{0}},g_{m}) $ of mass $ \mathfrak{m}=2r_{0}(k-1) $ with $ r_{0}=\frac{\mathfrak{c}_{\Sigma}}{k} $.
\end{theorem}
\begin{remark}
When $ M $ is nonparbolic and its Ricci curvature is bounded from below, the Cheng-Yau gradient estimate \cite{Cheng-Yau} implies $ |\nabla u|\to 0 $ at infinity.
\end{remark}

\begin{remark}
When $ k=1 $ and $ M $ is boundaryless, the argument in Theorem~\ref{Th1.1} can be used to refind the results in Colding-Minicozzi \cite{Cold4}. When $ k=2 $ and $ M $ has a compact, connected minimal boundary, Theorem~\ref{Th1.1} recovers the results in Oronzio \cite{Oro2}.
\end{remark}
Using the gradient estimates in Theorem~\ref{Th1.1}, we can derive a Cohn-Vossen type bound for the average of the scalar curvature. This generalizes the results in Colding-Minicozzi \cite{Cold4}.
\begin{corollary}\label{co1.1}
Under the assumptions in Theorem~\ref{Th1.1}, for any $ t\in [\frac{\mathfrak{c}_{\Sigma}}{k}, +\infty)  $, 
$$
t^{-1}\int_{\frac{\mathfrak{c}_{\Sigma}}{k}}^{t}\int_{\Sigma_{\tau}}R_{g}\leq 48\pi.
$$
\end{corollary}
This paper is organized as follows. In Section 2, we present some preliminary results. In Section 3, we first derive two differential inequalities and then establish a quadratic growth of the gradient integral estimates. In Section 4, we prove the main results.
\section{\bf Preliminaries}
Let $ (M,g) $ be a 3-dimensional, complete, noncompact and one-ended Riemannian manifold with a compact and connected boundary $ \Sigma $ and with nonnegative scalar curvature.  Let $ u $ be the solution to (\ref{u}) and let $ \Sigma_{t}=\{x\in M| u(x)=f(t)\} $ be the level set of $ u $. Along the level set $ \Sigma_{t} $, we define 
$$
\mathcal{A}_{k}(t)=\frac{t^{2}}{\mathfrak{c}_{\Sigma}^{2}}\Big(1+\frac{k-1}{kt}\mathfrak{c}_{\Sigma}\Big)^{4}\int_{\Sigma_{t}}|\nabla u|^{2},
$$

$$
\mathcal{G}_{k}(t)=-\frac{4\pi\mathfrak{c}_{\Sigma}^{2}}{k^{2}t}+\frac{t}{k^{2}}\Big(1+\frac{k-1}{kt}\mathfrak{c}_{\Sigma} \Big)^{4}\int_{\Sigma_{t}}|\nabla u|^{2},
$$

$$
\mathcal{F}_{k}(t)=4\pi t+\frac{t^{3}}{\mathfrak{c}_{\Sigma}^{2}}\Big(1+\frac{k-1}{kt}\mathfrak{c}_{\Sigma}\Big)^{3}\Big(1-\frac{3(k-1)}{kt}\mathfrak{c}_{\Sigma}\Big)\int_{\Sigma_{t}}|\nabla u|^{2}-\frac{t^{2}}{\mathfrak{c}_{\Sigma}}\Big(1+\frac{k-1}{kt}\mathfrak{c}_{\Sigma}\Big)^{2}\int_{\Sigma_{t}}|\nabla u|H.
$$
Set
$$
A_{0}=-\mathcal{G}_{k}(\frac{\mathfrak{c}_{\Sigma}}{k})=\mathfrak{c}_{\Sigma}\Big(\frac{4\pi}{k}-k\int_{\Sigma}|\nabla u|^{2}\Big),
$$

$$
\mathcal{F}_{k}(\frac{\mathfrak{c}_{\Sigma}}{k})=\mathfrak{c}_{\Sigma}\Big(\frac{4\pi}{k}+(4-3k)\int_{\Sigma}|\nabla u|^{2}-\int_{\Sigma}|\nabla u|H \Big).
$$
For $ \alpha\in (-\frac{1}{k\mathfrak{c}_{\Sigma}}, \frac{1}{k\mathfrak{c}_{\Sigma}} ) $, we denote 
\begin{equation}\label{A1}
A_{1}=(1-k\alpha \mathfrak{c}_{\Sigma})A_{0},
\end{equation}

$$
\widetilde{\mathcal{A}_{k}}(t)=\mathcal{A}_{k}(t)+\frac{1}{2t}A_{1}.
$$
We introduce the $ (0,2) $-symmetric tensors of $ M $:
$$
T=\nabla^{2}u+\frac{3(2\omega-k)}{(1-u)\omega}du\otimes du,\,\,\,B=T-\frac{tr(T)}{3}g,
$$
where $ \omega(u)=1+(k-1)u $. We define
$$
B_{1}(t)=\int_{\Sigma_{t}}\frac{|B|^{2}}{|\nabla u|^{2}},\,\,\,R_{1}(t)=\int_{\Sigma_{t}}R_{g}.
$$
Observe that 
$$
\frac{|B|^{2}}{2|\nabla u|^{2}}=\frac{|\nabla^{\Sigma_{t}}|\nabla u| |^{2}}{|\nabla u|^{2}}+\frac{1}{2}|\overset{\circ}{\mathbb{II}}|^{2}+\frac{3}{4}\Big(H-\frac{2(2\omega-k)}{(1-u)\omega}|\nabla u|\Big)^{2}.
$$
The monotone quantitity $ \mathcal{F}_{k}(t) $ was investigated in \cite{Ag2,Ag3} for p-harmonic function with $ p\in (1,3) $ and $ k=1 $. In \cite{Oro1}, Oronzio considered $ \mathcal{F}_{k}(t) $ for harmonic functions with $ k=2 $. In \cite{Chen}, Chen-Fu considered $ \mathcal{F}_{k}(t) $ for harmonic functions with $ k>0 $. In \cite{Xia}, Xia-Yin-Zhou derived $ \mathcal{F}_{k}(t) $ for p-harmonic functions with $ p\in (1,3) $ and $ k\in (0,2] $. Hence, we have the following result.
\begin{proposition}\cite{Oro1,Chen,Xia}\label{F mono}
Let $ (M,g) $ be a complete, orientable, one-ended 3-manifold with compact, connnected boundary $ \Sigma $. Suppose that $ H_{2}(M,\Sigma)=0 $. Let $ u $ be the solution of (\ref{u}). If $ R_{g}\geq 0 $ and $ k>0 $, then the following quantity 
\begin{equation}\label{F}
\mathcal{F}_{k}(t)=4\pi t+\frac{t^{3}}{\mathfrak{c}_{\Sigma}^{2}}\Big(1+\frac{k-1}{kt}\mathfrak{c}_{\Sigma}\Big)^{3}\Big(1-\frac{3(k-1)}{kt}\mathfrak{c}_{\Sigma}\Big)\int_{\Sigma_{t}}|\nabla u|^{2}-\frac{t^{2}}{\mathfrak{c}_{\Sigma}}\Big(1+\frac{k-1}{kt}\mathfrak{c}_{\Sigma}\Big)^{2}\int_{\Sigma_{t}}|\nabla u|H
\end{equation}
is monotone nondecreasing in $ t $. Futhermore, at a regular value $ t $, we have 
\begin{equation}\label{Fk}
\begin{aligned}
\mathcal{F}_{k}^{\prime}(t)&=4\pi-2\pi\chi(\Sigma_{t})+\frac{1}{2}\int_{\Sigma_{t}}\Big(\frac{|B|^{2}}{|\nabla u|^{2}}+R_{g} \Big).
\end{aligned}
\end{equation}
Moreover, if $ \mathcal{F}_{k}(t)\equiv0 $ for all  $ t\in [\frac{\mathfrak{c}_{\Sigma}}{k}, +\infty) $, then $ (M,g) $ is isometric to $ (M_{m,r_{0}},g_{m}) $ of mass $ \mathfrak{m}=2r_{0}(k-1) $ with $ r_{0}=\frac{\mathfrak{c}_{\Sigma}}{k} $.
\end{proposition}
\section{\bf Quadratic growth}
In this section, following an analogus argument in \cite{Cold4} and \cite{Oro2}, we first derive two differential inequalities and then establish a quadratic growth of $ \mathcal{A}_{k}(t) $ under the boundary mean curvature condition (\ref{H condition0 0}). This boundary mean curvature condition is inspired by the work of \cite{Chen}, \cite{Oro1} and \cite{Oro2}, and it is essential to obtain the quadratic growth of $ \mathcal{A}_{k}(t) $. 

Under the boundary mean curvature condition (\ref{H condition0 0}), we claim that 
\begin{equation}\label{boundary condition}
\mathcal{F}_{k}(\frac{\mathfrak{c}_{\Sigma}}{k})\geq (1-k\alpha \mathfrak{c}_{\Sigma})A_{0}=A_{1}.
\end{equation}
In fact, we have 
$$
\int_{\Sigma}H|\nabla u|\leq 4\pi \alpha \mathfrak{c}_{\Sigma}-k^{2}\alpha \mathfrak{c}_{\Sigma}\int_{\Sigma}|\nabla u|^{2}.
$$
Hence we get 
\begin{equation}\label{boundary con}
\begin{aligned}
\mathcal{F}_{k}(\frac{\mathfrak{c}_{\Sigma}}{k})&=\mathfrak{c}_{\Sigma}\Big(\frac{4\pi}{k}+(4-3k)\int_{\Sigma}|\nabla u|^{2}-\int_{\Sigma}|\nabla u|H \Big)\\
&\geq \mathfrak{c}_{\Sigma}\Big(\frac{4\pi}{k^{2}}(k-k^{2}\alpha \mathfrak{c}_{\Sigma})-(3k-4-k^{2}\alpha \mathfrak{c}_{\Sigma})\int_{\Sigma}|\nabla u|^{2}\Big)\\
&\geq \mathfrak{c}_{\Sigma}\Big(\frac{4\pi}{k^{2}}(k-k^{2}\alpha \mathfrak{c}_{\Sigma})-(k-k^{2}\alpha \mathfrak{c}_{\Sigma})\int_{\Sigma}|\nabla u|^{2} \Big)\\
&=\mathfrak{c}_{\Sigma}(1-k\alpha \mathfrak{c}_{\Sigma})\Big(\frac{4\pi}{k}-k\int_{\Sigma}|\nabla u|^{2} \Big)\\
&=(1-k\alpha \mathfrak{c}_{\Sigma})A_{0}=A_{1},
\end{aligned}
\end{equation}
where we have used the fact that $ 3k-4\leq k $ provided $ k\in (0,2] $.

Appling the estimate (\ref{boundary condition}), one has the following differential inequalities.
\begin{proposition}\label{diff ineq}
Let $ (M,g) $ be a complete, orientable, one-ended 3-manifold with compact, connnected boundary $ \Sigma $ and with nonnegative scalar curvature. Suppose that $ H_{2}(M,\Sigma)=0 $. Let $ u $ be the solution of (\ref{u}) and let $ k\in (0,2] $. If there exists 
\begin{equation}\label{a condi}
\alpha\in (-\frac{1}{k\mathfrak{c}_{\Sigma}},\frac{1}{k\mathfrak{c}_{\Sigma}}) \end{equation}
such that 
\begin{equation}\label{H condition0}
H\leq \alpha(1-k^{2}\mathfrak{c}_{\Sigma}|\nabla u|)\,\,\,\text{on }\Sigma.
\end{equation}
Then we have the following differential inequalities.
\\
$ (i) $ If $ A_{1}\geq 0 $, then 
$$
t\mathcal{A}_{k}^{\prime}(t)\geq \mathcal{A}_{k}(t)-4\pi+\frac{1}{2t}\int_{\frac{\mathfrak{c}_{\Sigma}}{k}}^{t}(B_{1}(s)+R_{1}(s))ds
$$
holds everywhere.
\\
$ (ii) $ For all $ t\in [\frac{\mathfrak{c}_{\Sigma}}{k}, +\infty) $, 
$$
t\widetilde{\mathcal{A}_{k}}^{\prime}(t)\geq \widetilde{\mathcal{A}_{k}}(t)-4\pi+\frac{1}{2t}\int_{\frac{\mathfrak{c}_{\Sigma}}{k}}^{t}(B_{1}(s)+R_{1}(s))ds
$$
holds everywhere.
\end{proposition}
\begin{proof}
At a regular value $ t\in \mathcal{T} $, we compute
$$
\mathcal{A}_{k}^{\prime}(t)=\frac{2t}{\mathfrak{c}_{\Sigma}^{2}}\Big(1+\frac{k-1}{kt}\mathfrak{c}_{\Sigma}\Big)^{3}\Big(1-\frac{k-1}{kt}\mathfrak{c}_{\Sigma}\Big)\int_{\Sigma}|\nabla u|^{2}-\frac{1}{\mathfrak{c}_{\Sigma}}\Big(1+\frac{k-1}{kt}\mathfrak{c}_{\Sigma}\Big)^{2}\int_{\Sigma_{t}}|\nabla u|H.
$$
Then we get 
\begin{equation}\label{t Ak}
\begin{aligned}
t\mathcal{A}_{k}^{\prime}(t)-\mathcal{A}_{k}(t)=&\frac{t^{2}}{\mathfrak{c}_{\Sigma}^{2}}\Big(1+\frac{k-1}{kt}\mathfrak{c}_{\Sigma}\Big)^{3}\Big(1-\frac{3(k-1)}{kt}\mathfrak{c}_{\Sigma}\Big)\int_{\Sigma}|\nabla u|^{2}\\
&-\frac{t}{\mathfrak{c}_{\Sigma}}\Big(1+\frac{k-1}{kt}\mathfrak{c}_{\Sigma}\Big)^{2}\int_{\Sigma_{t}}|\nabla u|H\\
=&\frac{1}{t}\mathcal{F}_{k}(t)-4\pi\\
=&\frac{1}{t}\Big(\mathcal{F}_{k}(\frac{\mathfrak{c}_{\Sigma}}{k})+\int_{\frac{\mathfrak{c}_{\Sigma}}{k}}^{t} \mathcal{F}_{k}^{\prime}(s)ds\Big)-4\pi\\
\geq &\frac{1}{t}A_{1}+\frac{1}{t}\int_{\frac{\mathfrak{c}_{\Sigma}}{k}}^{t} \mathcal{F}_{k}^{\prime}(s)ds-4\pi,
\end{aligned}
\end{equation}
where we have used (\ref{boundary condition}).
On the other hand, by Guass-Bonnet theorem and (\ref{Fk}) , we have
\begin{equation}
\begin{aligned}\label{F ineq}
\mathcal{F}_{k}^{\prime}(t)&=4\pi-2\pi\chi(\Sigma_{t})+\frac{1}{2}\int_{\Sigma_{t}}\Big(\frac{|B|^{2}}{|\nabla u|^{2}}+R_{g} \Big)\\
&\geq \frac{1}{2}(B_{1}(t)+R_{1}(t)).
\end{aligned}
\end{equation}
Thus we find 
\begin{equation}\label{diff Ak}
t\mathcal{A}_{k}^{\prime}(t)\geq \mathcal{A}_{k}(t)-4\pi
+\frac{A_{1}}{t}+\frac{1}{2t}\int_{\frac{\mathfrak{c}_{\Sigma}}{k}}^{t}(B_{1}(s)+R_{1}(s))ds.
\end{equation}
By Sard's theorem and continuity of $ \mathcal{A}_{k}(t) $, (\ref{diff Ak}) holds for every $ t\in [\frac{\mathfrak{c}_{\Sigma}}{k}, +\infty) $. If $ A_{1}\geq 0 $, we obtain $ (i) $.

To prove $ (ii) $, a simple calculation shows that
$$
\begin{aligned}
t\widetilde{\mathcal{A}_{k}}^{\prime}(t)-\widetilde{\mathcal{A}_{k}}(t)&=t\mathcal{A}_{k}^{\prime}(t)-\mathcal{A}_{k}(t)-\frac{1}{t}A_{1}\\
&\geq -4\pi+\frac{1}{2t}\int_{\frac{\mathfrak{c}_{\Sigma}}{k}}^{t}(B_{1}(s)+R_{1}(s))ds.
\end{aligned}
$$
This proves $ (ii) $.
\end{proof}
Using Proposition~\ref{diff ineq}, we now establish the quadratic growth of $ \mathcal{A}_{k}(t) $. 
\begin{proposition}\label{quadratic growth}
Under the assumptions in Proposition~\ref{diff ineq}, we have the following
\\
$ (i) $ If $ A_{1}\geq 0 $ and $ \mathcal{A}_{k}(t_{0})>4\pi $ for some $ t_{0}\in [\frac{\mathfrak{c}_{\Sigma}}{k},+\infty ) $, then there exists a positive constant $ C>0 $ such that, for every $ t\geq t_{0} $,
$$
\mathcal{A}_{k}(t)\geq Ct^{2}.
$$
\\
$ (ii) $ If $ A_{1}<0 $ and $ \mathcal{A}_{k}(t_{0})>4\pi-\frac{A_{1}}{t_{0}} $ for some $ t_{0}\in [\frac{\mathfrak{c}_{\Sigma}}{k},+\infty ) $, then there exists a positive constant $ C>0 $ such that, for every $ t\geq t_{0} $,
$$
\mathcal{A}_{k}(t)\geq Ct^{2}.
$$
\end{proposition}
\begin{remark}
When $ k=1 $ and $ M $ is boundaryless, it suffices to consider the case $ (i) $ in Proposition~\ref{diff ineq} and Proposition~\ref{quadratic growth}. When $ k=2 $ and $ M $ has compact, connected minimal boundary, Proposition~\ref{diff ineq} and Proposition~\ref{quadratic growth} recover the results in Oronzio \cite{Oro2}.
\end{remark}
\begin{proof}
First we choose $ \delta>0 $ such that 
\begin{equation}\label{Ak t0}
\mathcal{A}_{k}(t_{0})>4\pi+\delta,\,\,\,\text{if }A_{1}\geq 0,
\end{equation}
and 
\begin{equation}\label{wAk t0}
\widetilde{\mathcal{A}_{k}}(t_{0})>4\pi+\delta,\,\,\,\text{if }A_{1}<0.
\end{equation}
In fact, if $ A_{1}\geq0 $ and $ \mathcal{A}_{k}(t_{0})>4\pi $, we can choose $ \delta>0 $ small enough such that (\ref{Ak t0}) holds. If $ A_{1}<0 $ and $ \mathcal{A}_{k}(t_{0})>4\pi-\frac{A_{1}}{t_{0}} $, we deduce that $ \widetilde{\mathcal{A}_{k}}(t_{0})>4\pi $, so we can again choose $ \delta>0 $ small enough such that (\ref{wAk t0}) holds.

In both cases, we claim that there exists $ \delta>0 $, such that for every $ t\geq t_{0} $,
\begin{equation}\label{key estimate}
\mathcal{A}_{k}^{\prime}(t)\geq 0,\,\,\,\,\mathcal{A}_{k}(t)\geq 4\pi+\delta\frac{t}{t_{0}}.
\end{equation}
Indeed, let $ \delta>0 $ be chosen as in (\ref{Ak t0}) and (\ref{wAk t0}). For case $ (i) $, by Proposition~\ref{diff ineq} $ (i) $, we obtain that 
\begin{equation}\label{diff ineq case1}
t\mathcal{A}_{k}^{\prime}(t)\geq \mathcal{A}_{k}(t)-4\pi.
\end{equation}
In this case, $ \mathcal{A}_{k}(t_{0})>4\pi $, so that $ \mathcal{A}_{k}^{\prime}(t_{0})>0 $. This implies that 
\begin{equation}\label{A A'}
\mathcal{A}_{k}(t)>4\pi,\,\,\,\,\mathcal{A}_{k}^{\prime}(t)\geq0
\end{equation}
for all $ t\geq t_{0} $. Indeed, we consider the function $ h(t)=\frac{1}{t}(\mathcal{A}_{k}(t)-4\pi) $, then (\ref{diff ineq case1}) gives 
$$
h^{\prime}(t)=\frac{1}{t^{2}}\Big(t\mathcal{A}_{k}^{\prime}(t)- (\mathcal{A}_{k}(t)-4\pi)\Big)\geq 0.
$$
Hence $ h(t)\geq h(t_{0})>0 $, that is $ \mathcal{A}_{k}(t)>4\pi $, and so $ \mathcal{A}_{k}^{\prime}(t)\geq0 $.
Now (\ref{diff ineq case1}) can be rewritten as 
$$
\frac{\mathcal{A}_{k}^{\prime}(t)}{\mathcal{A}_{k}(t)-4\pi}\geq \frac{1}{t}
$$
for all $ t\geq t_{0} $. Integrating from $ t_{0} $ to $ t $ gives 
$$
\mathcal{A}_{k}(t)-4\pi\geq (\mathcal{A}_{k}(t_{0})-4\pi)\frac{t}{t_{0}}>\delta\frac{t}{t_{0}}.
$$
Hence (\ref{key estimate}) holds.

For case $ (ii) $, by (\ref{diff Ak}), we have 
\begin{equation}\label{diff ineq case2}
t\mathcal{A}_{k}^{\prime}(t)\geq \mathcal{A}_{k}(t)-(4\pi
-\frac{A_{1}}{t}).
\end{equation}
In this case, $ \mathcal{A}_{k}(t_{0})>4\pi-\frac{A_{1}}{t_{0}} $, so that $ \mathcal{A}_{k}^{\prime}(t_{0})>0 $. This shows that 
$$
\mathcal{A}_{k}(t)>4\pi-\frac{A_{1}}{t}
$$
for all $ t\geq t_{0} $, which implies that $ \mathcal{A}_{k}^{\prime}(t)\geq0 $ for all $ t\geq t_{0} $. Indeed, we consider the function $ g(t)=\frac{1}{t}(\mathcal{A}_{k}(t)-4\pi+\frac{A_{1}}{t}) $, then (\ref{diff ineq case2}) gives 
$$
g^{\prime}(t)=\frac{1}{t^{2}}\Big(t\mathcal{A}_{k}^{\prime}(t)-(\mathcal{A}_{k}(t)-4\pi)-\frac{2A_{1}}{t} \Big)\geq \frac{1}{t^{2}}\Big(t\mathcal{A}_{k}^{\prime}(t)-(\mathcal{A}_{k}(t)-4\pi)-\frac{A_{1}}{t} \Big)\geq 0.
$$
due to $ A_{1}<0 $. Hence $ g(t)\geq g(t_{0})>0 $, that is $ \mathcal{A}_{k}(t)>4\pi-\frac{A_{1}}{t} $. On the other hand, by Proposition~\ref{diff ineq} $ (ii) $, we have 
$$
t\widetilde{\mathcal{A}_{k}}^{\prime}(t)\geq \widetilde{\mathcal{A}_{k}}(t)-4\pi.
$$
Since $ \widetilde{\mathcal{A}_{k}}(t_{0})>4\pi $, we deduce that $ \widetilde{\mathcal{A}_{k}}^{\prime}(t_{0})>0 $, and so 
$$
\widetilde{\mathcal{A}_{k}}(t)>4\pi,\,\,\,\widetilde{\mathcal{A}_{k}}^{\prime}(t)\geq 0,
$$
for all $ t\geq t_{0} $. Arguing as in case $ (i) $ and using (\ref{wAk t0}), one has 
$$
\widetilde{\mathcal{A}_{k}}(t)\geq 4\pi+\delta\frac{t}{t_{0}}
$$
for all $ t\geq t_{0} $.
Since $ A_{1}<0 $, we have
$$ \widetilde{\mathcal{A}_{k}}(t)=\mathcal{A}_{k}(t)+\frac{1}{2t}A_{1}<\mathcal{A}_{k}(t).
$$ 
Therefore, we conclude that 
$$
\mathcal{A}_{k}(t)\geq 4\pi+\delta\frac{t}{t_{0}}
$$
for all $ t\geq t_{0} $. Hence (\ref{key estimate}) holds.

As in \cite{Cold4,Oro2}, on the interval $ [\frac{\mathfrak{c}_{\Sigma}}{k}, +\infty) $, we define the function 
$$
a(t)=t\frac{\mathcal{A}_{k}^{\prime}(t)}{\mathcal{A}_{k}(t)}.
$$
From (\ref{key estimate}), we see that $ a(t)\geq0 $ for all $ t\geq t_{0} $. By (\ref{t Ak}), for every $ t\in \mathcal{T} $, one has
$$
t\mathcal{A}_{k}^{\prime}(t)-\mathcal{A}_{k}(t)=\frac{1}{t}\mathcal{F}_{k}(t)-4\pi.
$$ 
This leads to 
$$
\begin{aligned}
\mathcal{A}_{k}^{\prime}(t)&=\frac{\mathcal{F}_{k}(t)}{t^{2}}+\frac{1}{t}(\mathcal{A}_{k}(t)-4\pi),\\
t\mathcal{A}_{k}^{\prime\prime}(t)&=[t\mathcal{A}_{k}^{\prime}(t)-\mathcal{A}_{k}(t)]^{\prime}=-\frac{\mathcal{F}_{k}(t)}{t^{2}}+\frac{\mathcal{F}_{k}^{\prime}(t)}{t}.
\end{aligned}
$$
Hence one has 
$$
a^{\prime}(t)=\frac{1}{t}\Big(1-\frac{4\pi}{\mathcal{A}_{k}(t)}+\frac{\mathcal{F}_{k}^{\prime}(t)}{\mathcal{A}_{k}(t)}-a^{2}(t) \Big).
$$
For $ \nu=\frac{\nabla u}{|\nabla u|} $, observe that 
$$
\begin{aligned}
B(\nu,\nu)&=-H|\nabla u|+\frac{2(2\omega(u)-k)}{(1-u)\omega(u)}|\nabla u|^{2}\\
&=-H|\nabla u|+\frac{2}{\mathfrak{c}_{\Sigma}}\Big(1+\frac{k-1}{kt}\mathfrak{c}_{\Sigma}\Big)\Big(1-\frac{k-1}{kt}\mathfrak{c}_{\Sigma} \Big)|\nabla u|^{2}.
\end{aligned}
$$
Then 
$$
t\mathcal{A}_{k}^{\prime}(t)=\frac{t}{\mathfrak{c}_{\Sigma}}\Big(1+\frac{k-1}{kt}\mathfrak{c}_{\Sigma}\Big)^{2}\int_{\Sigma_{t}}B(\nu,\nu).
$$
Using H\"older's inequality and refined Kato's inequality in \cite{Cold4} , we get 
$$
\begin{aligned}
t^{2}(\mathcal{A}_{k}^{\prime}(t))^{2}&=\frac{t^{2}}{\mathfrak{c}_{\Sigma}^{2}}\Big(1+\frac{k-1}{kt}\mathfrak{c}_{\Sigma} \Big)^{4}\Big(\int_{\Sigma_{t}}B(\nu,\nu)\Big)^{2}\\
&\leq \frac{t^{2}}{\mathfrak{c}_{\Sigma}^{2}}\Big(1+\frac{k-1}{kt}\mathfrak{c}_{\Sigma} \Big)^{4}\Big(\int_{\Sigma_{t}}|\nabla u|^{2} \Big)\Big(\int_{\Sigma_{t}}\frac{B(\nu,\nu)^{2}}{|\nabla u|^{2}}\Big)\\
&=\mathcal{A}_{k}(t)\int_{\Sigma_{t}}\frac{B(\nu,\nu)^{2}}{|\nabla u|^{2}}\\
&\leq \frac{2}{3}\mathcal{A}_{k}(t)B_{1}(t).
\end{aligned}
$$
On the other hand, by (\ref{F ineq}) and $ R_{g}\geq 0 $, we have 
$$
\mathcal{F}_{k}^{\prime}(t)\geq \frac{1}{2}B_{1}(t).
$$
Thus one has 
$$
\frac{\mathcal{F}_{k}^{\prime}(t)}{\mathcal{A}_{k}(t)}\geq \frac{B_{1}(t)}{2\mathcal{A}_{k}(t)}\geq \frac{3}{4}a(t)^{2}.
$$
Therefore, 
\begin{equation}\label{a ineq}
a^{\prime}(t)\geq \frac{1}{t}\Big(1-\frac{4\pi}{\mathcal{A}_{k}(t)}-\frac{1}{4}a(t)^{2} \Big)
\end{equation}
for each $ t\in \mathcal{T} $.
Now, we show the existence of a positive constant $ C=C(t_{0},\delta) $ such that 
$$
a(t)\geq 2-\frac{C}{\sqrt{t}}
$$
for all $ t\geq t_{0} $. If $ a(t)\geq 2 $ for all $ t\geq t_{0} $, the claim is trivially true. Therefore, we suppose that there exists $ \bar{t}_{0}\geq t_{0} $ such that 
$$
a(\bar{t}_{0})<2.
$$
Without loss of generality, we may assume that $ \bar{t}_{0}>t_{0} $ by the continuity of $ a(t) $ and that $ \bar{t}_{0}\in \mathcal{T} $ by Sard's theorem. 

Let $ t_{2}>t_{0} $ be a generic regular value such that $ a(t_{2})<2 $. Denote 
$$
t_{1}=\inf\limits_{t\in [t_{0},t_{2}]}\{t\,\,|\,\,a(s)<2,\,\,\forall s\in (t,t_{2}] \}.
$$ 
Either $ t_{1}=t_{0} $ or $ t_{1}>t_{0} $. We will consider two cases.

Suppose first that $ t_{1}=t_{0} $. In this case, one has 
\begin{equation}\label{a(t) 0 2}
0\leq a(t)\leq 2,\,\,\,\forall t\in [t_{0}, t_{2}].
\end{equation}
Hence, it follows from (\ref{a ineq}) that 
$$
a^{\prime}(t)\geq \frac{1}{t}\Big(1-\frac{4\pi}{\mathcal{A}_{k}(t)}-\frac{1}{2}a(t) \Big),\,\,\,\,\,\forall t\in [t_{0},t_{2}].
$$
This implies that 
\begin{equation}\label{a ineq1}
(\sqrt{t}a(t))^{\prime}\geq \frac{1}{\sqrt{t}}\Big(1-\frac{4\pi}{\mathcal{A}_{k}(t)}\Big)
\end{equation}
for every $ t\in \mathcal{T}\cap [t_{0},t_{2}] $. Observe that the function $ \sqrt{t}a(t) $ is locally absolute continuous in $ [\frac{\mathfrak{c}_{\Sigma}}{k},+\infty)$. Thus, integrating (\ref{a ineq1}) from $ t_{0} $ to $ t_{2} $ gives 
$$
a(t_{2})\geq 2+\sqrt{\frac{t_{0}}{t_{2}}}a(t_{0})-2\sqrt{\frac{t_{0}}{t_{2}}}-\frac{4\pi}{\sqrt{t_{2}}}\int_{t_{0}}^{t_{2}}\frac{ds}{\sqrt{s}\mathcal{A}_{k}(s)}.
$$ 
Then, it follows from (\ref{key estimate}) and (\ref{a(t) 0 2}) that 
$$
\begin{aligned}
a(t_{2})&\geq 2+\sqrt{\frac{t_{0}}{t_{2}}}a(t_{0})-2\sqrt{\frac{t_{0}}{t_{2}}}-\frac{4\pi t_{0}}{\delta \sqrt{t_{2}}}\int_{t_{0}}^{t_{2}}\frac{1}{s^{\frac{3}{2}}}ds\\
&=2+\sqrt{\frac{t_{0}}{t_{2}}}a(t_{0})-2\sqrt{\frac{t_{0}}{t_{2}}}-\frac{8\pi t_{0}}{\delta \sqrt{t_{2}}}(t_{0}^{-\frac{1}{2}}-t_{2}^{-\frac{1}{2}})\\
&\geq 2+\sqrt{\frac{t_{0}}{t_{2}}}a(t_{0})-2\sqrt{\frac{t_{0}}{t_{2}}}-\frac{8\pi}{\delta}\sqrt{\frac{t_{0}}{t_{2}}}\\
&\geq 2-2(1+\frac{4\pi}{\delta})\sqrt{\frac{t_{0}}{t_{2}}}.
\end{aligned}
$$
In the remaining case, we have  $ t_{1}>t_{0} $. Thus one deduces that $ a(t_{1})=2 $ and that $ 0\leq a(t)\leq 2 $, for any $ t\in [t_{1},t_{2}] $. Hence, proceeding as above, we get 
$$
\begin{aligned}
a(t_{2})&\geq 2+\sqrt{\frac{t_{1}}{t_{2}}}a(t_{1})-2\sqrt{\frac{t_{1}}{t_{2}}}-\frac{4\pi}{\sqrt{t_{2}}}\int_{t_{1}}^{t_{2}}\frac{ds}{\sqrt{s}\mathcal{A}_{k}(s)}\\
&=2-\frac{4\pi}{\sqrt{t_{2}}}\int_{t_{1}}^{t_{2}}\frac{ds}{\sqrt{s}\mathcal{A}_{k}(s)}\\
&\geq 2-\frac{4\pi}{\sqrt{t_{2}}}\int_{t_{0}}^{t_{2}}\frac{ds}{\sqrt{s}\mathcal{A}_{k}(s)}\\
&\geq 2-\frac{8\pi}{\delta}\sqrt{\frac{t_{0}}{t_{2}}}.
\end{aligned}
$$
Therefore, choosing $ C=C(t_{0},\delta)=2(1+\frac{4\pi}{\delta})\sqrt{t_{0}} $, one has 
\begin{equation}\label{a regul}
a(t)\geq 2-\frac{C}{\sqrt{t}},
\end{equation}
for every $ t\in \mathcal{T}\cap [t_{0},+\infty) $. The continuity of the function $ a(t) $ implies that (\ref{a regul}) holds everywhere in the interval $ [t_{0},+\infty) $. Then, by the definition of $ a(t) $, we have 
$$
(\log(\mathcal{A}_{k}(t)))^{\prime}\geq 2t^{-1}-Ct^{-\frac{3}{2}}
$$ 
for all $ t\in [t_{0},+\infty) $. Integrating this inequality from $ t_{0} $ to $ t $, we obtain that 
$$
\mathcal{A}_{k}(t)\geq t_{0}^{-2}\mathcal{A}_{k}(t_{0})e^{-\frac{2C}{\sqrt{t_{0}}}}t^{2},
$$
which is the required quadratic growth of $ \mathcal{A}_{k}(t) $.
\end{proof}
\section{Proof of main results}
In this section, we aim to prove Theorem~\ref{Th1.1} and Corollary~\ref{co1.1}. To prove Theorem~\ref{Th1.1}, we first use Proposition~\ref{quadratic growth} to present a contradiction argument to obtain $ (a) $ and $ (b) $ in Theorem~\ref{Th1.1}. Then, we use the monotonicity and rigidity of $ \mathcal{G}_{k}(t) $ and $ \mathcal{F}_{k}(t) $ to prove the remaining part of Theorem~\ref{Th1.1}. Finally, we apply the gradient estimate in Theorem~\ref{Th1.1} to prove Corollary~\ref{co1.1}.

\textbf{Proof of Theorem 1.1}. First we assume by contradiction that $ A_{0}<0 $. Then, by (\ref{a condi 0}) and (\ref{A1}), we get that $ A_{1}<0 $. We will consider two cases. Either 
$$
\mathcal{A}_{k}(t)\leq 4\pi-\frac{A_{1}}{t}
$$
for all $ t\in [\frac{\mathfrak{c}_{\Sigma}}{k},+\infty) $ or there exists some $ t_{0}\in [\frac{\mathfrak{c}_{\Sigma}}{k},+\infty) $ such that 
$$
\mathcal{A}_{k}(t_{0})> 4\pi-\frac{A_{1}}{t_{0}}.
$$
We claim that the last case is not possible. Indeed, if this holds, Proposition~\ref{quadratic growth} implies that there exists a positive constant $ C>0 $ such that 
$$
\mathcal{A}_{k}(t)\geq Ct^{2}
$$ 
for all $ t\geq t_{0} $, which can be rewritten as 
$$
\int_{\Sigma_{t}}|\nabla u|^{2}\geq C\Big(1+\frac{k-1}{kt}\mathfrak{c}_{\Sigma}\Big)^{-4}\mathfrak{c}_{\Sigma}^{2}\geq k^{-4}C\mathfrak{c}_{\Sigma}^{2}.
$$
This is a contradiction, since $ |\nabla u|\to 0 $ at infinity and $ \int_{\Sigma_{t}}|\nabla u|=4\pi \mathfrak{c}_{\Sigma} $ for all $ t\in [\frac{\mathfrak{c}_{\Sigma}}{k}, +\infty) $.
Therefore, we have 
$$
\mathcal{A}_{k}(t)\leq 4\pi-\frac{A_{1}}{t}
$$
for all $ t\in [\frac{\mathfrak{c}_{\Sigma}}{k}, +\infty) $. This implies that 
$$
\lim\limits_{t\to +\infty}t\Big(1+\frac{k-1}{kt}\mathfrak{c}_{\Sigma}\Big)^{4}\int_{\Sigma_{t}}|\nabla u|^{2}=0.
$$
Hence, we obtain that 
\begin{equation}\label{G lim}
\lim\limits_{t\to +\infty}\mathcal{G}_{k}(t)=\lim\limits_{t\to +\infty}\Big[-\frac{4\pi\mathfrak{c}_{\Sigma}^{2}}{k^{2}t}+\frac{t}{k^{2}}\Big(1+\frac{k-1}{kt}\mathfrak{c}_{\Sigma} \Big)^{4}\int_{\Sigma_{t}}|\nabla u|^{2}\Big]=0.
\end{equation}
For each $ t\in \mathcal{T} $, one has
\begin{equation}\label{G diff}
\mathcal{G}_{k}^{\prime}(t)=\frac{\mathfrak{c}_{\Sigma}^{2}}{k^{2}t^{3}}\mathcal{F}_{k}(t)\geq \frac{\mathfrak{c}_{\Sigma}^{2}}{k^{2}t^{3}}\mathcal{F}_{k}(\frac{\mathfrak{c}_{\Sigma}}{k})\geq \frac{\mathfrak{c}_{\Sigma}^{2}}{k^{2}t^{3}}A_{1}.
\end{equation}
Integrating (\ref{G diff}) from $ \frac{\mathfrak{c}_{\Sigma}}{k} $ to $ t $ gives 
$$
\mathcal{G}_{k}(t)-\mathcal{G}_{k}(\frac{\mathfrak{c}_{\Sigma}}{k})\geq \Big(\frac{1}{2}-\frac{\mathfrak{c}_{\Sigma}^{2}}{2k^{2}t^{2}}\Big)A_{1}.
$$
Letting $ t\to +\infty $ and using (\ref{G lim}), we get 
$$
A_{0}=-\mathcal{G}_{k}(\frac{\mathfrak{c}_{\Sigma}}{k})\geq \frac{1}{2}A_{1}=\frac{1}{2}(1-k\alpha\mathfrak{c}_{\Sigma})A_{0},
$$
thus
$$
\frac{1}{2}(1+k\alpha\mathfrak{c}_{\Sigma})A_{0}\geq 0.
$$
By (\ref{a condi 0}), we deduce that $ A_{0}\geq 0 $, which is again a contradiction.

In conclusion, we obtain that 
$$ 
A_{0}=-\mathcal{G}_{k}(\frac{\mathfrak{c}_{\Sigma}}{k})=\mathfrak{c}_{\Sigma}\Big(\frac{4\pi}{k}-k\int_{\Sigma_{t}}|\nabla u|^{2}\Big)\geq 0, 
$$ 
which implies the first part of the statement in $ (a) $. 

If $ \mathcal{A}_{k}(t_{0})>4\pi $ for some $ t_{0}\in [\frac{\mathfrak{c}_{\Sigma}}{k}, +\infty)  $, then use Proposition~\ref{quadratic growth} again to deduce that there exists a positive constant $ C>0 $ such that 
$$ 
\mathcal{A}_{k}(t)\geq Ct^{2} 
$$ 
for all $ t\geq t_{0} $. Arguing as before, we obtain a contradiction.
Therefore, we get 
\begin{equation}\label{A 4pi}
\mathcal{A}_{k}(t)\leq 4\pi
\end{equation}
for all $ t\in [\frac{\mathfrak{c}_{\Sigma}}{k}, +\infty) $. This implies that 
$$
t\Big(1+\frac{k-1}{kt}\mathfrak{c}_{\Sigma}\Big)^{4}\int_{\Sigma_{t}}|\nabla u|^{2}\leq \frac{4\pi \mathfrak{c}_{\Sigma}^{2}}{t},
$$
and so 
$$
\lim\limits_{t\to +\infty}\mathcal{G}_{k}(t)=0.
$$
Since $ A_{1}\geq0 $, by (\ref{G diff}), we conclude that $ \mathcal{G}_{k}(t) $ is monotone nondecreasing in $ t $. If $ A_{0}=0 $, then $ A_{1}=0 $ and 
$$
\lim\limits_{t\to+\infty}\mathcal{G}_{k}(t)=\mathcal{G}_{k}(\frac{\mathfrak{c}_{\Sigma}}{k})=-A_{0}=0, 
$$ 
which implies that 
$$
\mathcal{G}_{k}^{\prime}(t)=\frac{\mathfrak{c}_{\Sigma}^{2}}{k^{2}t^{3}}\mathcal{F}_{k}(t)\equiv0
$$
for all  $ t\in [\frac{\mathfrak{c}_{\Sigma}}{k}, +\infty) $. Now the rest part of $ (a) $ follows from Proposition~\ref{F mono}. 

To prove $ (b) $, the first part of the statement is contained in $ (\ref{A 4pi}) $. If $ \mathcal{A}_{k}(t_{0})=4\pi $ for some $ t_{0}\in [\frac{\mathfrak{c}_{\Sigma}}{k}, +\infty) $, then 
$$
\mathcal{G}_{k}(t_{0})=-\frac{4\pi\mathfrak{c}_{\Sigma}^{2}}{k^{2}t_{0}}+\frac{t_{0}}{k^{2}}\Big(1+\frac{k-1}{kt_{0}}\mathfrak{c}_{\Sigma} \Big)^{4}\int_{\Sigma_{t_{0}}}|\nabla u|^{2}=0,
$$
and so 
$$
\lim\limits_{t\to+\infty}\mathcal{G}_{k}(t)=\mathcal{G}_{k}(t_{0})=0.
$$
Arguing as before, we get 
\begin{equation}\label{G F 0}
\mathcal{G}_{k}(t)=\mathcal{F}_{k}(t)=\mathcal{F}_{k}^{\prime}(t)\equiv0
\end{equation}
for all $ t\geq t_{0} $. Since $ \mathcal{F}_{k}(t) $ is monotone nondecreasing in $ t $, by $ (a) $ and (\ref{G F 0}), we have 
$$
0\leq (1-k\alpha \mathfrak{c}_{\Sigma})A_{0}\leq \mathcal{F}_{k}(\frac{\mathfrak{c}_{\Sigma}}{k})\leq \lim\limits_{t\to \infty}\mathcal{F}_{k}(t)=0.
$$
Thus, $\mathcal{F}_{k}(t)\equiv0$ for all  $ t\in [\frac{\mathfrak{c}_{\Sigma}}{k}, +\infty) $. Now the second part of the statement in $ (b) $ follows from Proposition~\ref{F mono}.

To prove $ (c) $ and $ (d) $, observe that (\ref{A 4pi}) implies 
$$
\int_{\Sigma_{t}}|\nabla u|^{2}\leq \frac{4\pi \mathfrak{c}_{\Sigma}^{2}}{t^{2}\Big(1+\frac{k-1}{kt}\mathfrak{c}_{\Sigma}\Big)^{4}}.
$$
Recall that $ \Sigma_{t}^{*}=\Sigma_{t}\setminus\{|\nabla u|=0\} $, so we have 
$$
4\pi \mathfrak{c}_{\Sigma}=\int_{\Sigma_{t}}|\nabla u|=\int_{\Sigma_{t}^{*}}|\nabla u|\leq \Big(\int_{\Sigma_{t}^{*}}|\nabla u|^{2} \Big)^{\frac{1}{2}}|\Sigma_{t}^{*}|^{\frac{1}{2}}\leq \frac{(4\pi)^{\frac{1}{2}}\mathfrak{c}_{\Sigma}|\Sigma_{t}^{*}|^{\frac{1}{2}}}{t\Big(1+\frac{k-1}{kt}\mathfrak{c}_{\Sigma}\Big)^{2}}.
$$
Hence
$$
|\Sigma_{t}^{*}|\geq 4\pi t^{2}\Big(1+\frac{k-1}{kt}\mathfrak{c}_{\Sigma}\Big)^{4}.
$$
In particular, taking $ t=\frac{\mathfrak{c}_{\Sigma}}{k} $ and observing that $ \Sigma $ is regular, one has 
$$
|\Sigma|\geq 4\pi k^{2}\mathfrak{c}_{\Sigma}^{2}.
$$
Therefore, we have proved the first part of the statement of $ (c) $ and $ (d) $. Now the equality cases in  $ (c) $ and $ (d) $ follow from the equality cases in $ (b) $ and $ (a) $, respectively.

To prove $ (e) $, use $ (c) $ to obtain 
$$
\begin{aligned}
4\pi t^{2}\Big(1+\frac{k-1}{kt}\mathfrak{c}_{\Sigma}\Big)^{4}&\leq |\Sigma_{t}^{*}|=\int_{\Sigma_{t}^{*}}|\nabla u|^{\frac{2}{3}}\frac{1}{|\nabla u|^{\frac{2}{3}}}\\
&\leq \Big(\int_{\Sigma_{t}^{*}}|\nabla u|^{2} \Big)^{\frac{1}{3}}\Big(\int_{\Sigma_{t}^{*}}\frac{1}{|\nabla u|} \Big)^{\frac{2}{3}}\\
&\leq \Big[ \frac{4\pi \mathfrak{c}_{\Sigma}^{2}}{t^{2}\Big(1+\frac{k-1}{kt}\mathfrak{c}_{\Sigma}\Big)^{4}} \Big]^{\frac{1}{3}}\Big(\int_{\Sigma_{t}^{*}}\frac{1}{|\nabla u|} \Big)^{\frac{2}{3}}.
\end{aligned}
$$
Thus we get 
$$
\int_{\Sigma_{t}^{*}}\frac{1}{|\nabla u|}\geq 4\pi \mathfrak{c}_{\Sigma}^{-1}t^{4}\Big(1+\frac{k-1}{kt}\mathfrak{c}_{\Sigma}\Big)^{8}.
$$
By the coarea formula and the fact that $ \{|\nabla u|=0\} $ is negligible, one has 
$$
Vol\{u\leq f(t)\}=\int_{\frac{\mathfrak{c}_{\Sigma}}{k}}^{t}\frac{\mathfrak{c}_{\Sigma}}{s^{2}\Big(1+\frac{k-1}{ks}\mathfrak{c}_{\Sigma}\Big)^{2}}\int_{\Sigma_{s}}\frac{1}{|\nabla u|}ds\geq \int_{\frac{\mathfrak{c}_{\Sigma}}{k}}^{t}4\pi s^{2}\Big(1+\frac{k-1}{ks}\mathfrak{c}_{\Sigma}\Big)^{6}.
$$
The equality case follows from the equality case in $ (c) $.
$\hfill\square$ 

\textbf{Proof of Corollary 1.1}. From Theorem~\ref{Th1.1} $ (a) $, we know that $ A_{1}\geq0 $, it follows from (\ref{diff Ak}) that 
\begin{equation}\label{co Ak}
t\mathcal{A}_{k}^{\prime}(t)\geq \mathcal{A}_{k}(t)-4\pi
+\frac{1}{2t}\int_{\frac{\mathfrak{c}_{\Sigma}}{k}}^{t}(B_{1}(s)+R_{1}(s))ds.
\end{equation}
Observe that Theorem~\ref{Th1.1} $ (b) $ gives $ \mathcal{A}_{k}(t)\leq 4\pi $, so that
$$
\int_{t}^{2t}\mathcal{A}_{k}^{\prime}(s)ds=\mathcal{A}_{k}(2t)-\mathcal{A}_{k}(t)\leq 4\pi.
$$
Hence, we can take $ s\in [t,2t] $ with $ t\mathcal{A}_{k}^{\prime}(s)\leq 4\pi $. Using (\ref{co Ak}), one has 
$$
\begin{aligned}
8\pi&\geq 2t\mathcal{A}_{k}^{\prime}(s)\geq s\mathcal{A}_{k}^{\prime}(s)\geq \mathcal{A}_{k}(s)-4\pi
+\frac{1}{2s}\int_{\frac{\mathfrak{c}_{\Sigma}}{k}}^{s}(B_{1}(\tau)+R_{1}(\tau))d\tau\\
&\geq -4\pi+\frac{1}{2s}\int_{\frac{\mathfrak{c}_{\Sigma}}{k}}^{t}R_{1}(\tau)d\tau,
\end{aligned}
$$
that is
$$
\frac{1}{2s}\int_{\frac{\mathfrak{c}_{\Sigma}}{k}}^{t}R_{1}(\tau)d\tau\leq 12\pi.
$$
Therefore, we get 
$$
t^{-1}\int_{\frac{\mathfrak{c}_{\Sigma}}{k}}^{t}\int_{\Sigma_{\tau}}R_{g}= \Big(\frac{2s}{t}\Big)\frac{1}{2s}\int_{\frac{\mathfrak{c}_{\Sigma}}{k}}^{t}R_{1}(\tau)d\tau\leq 4(12\pi)=48\pi.
$$
$\hfill\square$ 

\end{document}